\documentclass[12pt]{amsart}

\usepackage{graphicx,amssymb,latexsym,amsfonts,txfonts,amsmath,amsthm}
\usepackage{pdfsync,color,tabularx,rotating}
\usepackage[all,cmtip]{xy}

\usepackage{tikz}

\advance\hoffset by -0.9 truecm

\newtheorem{theo}{Theorem}

\theoremstyle{remark}

\theoremstyle{remark}

\theoremstyle{remark}

\begin{document}

\title{A short proof of Greenberg's Theorem}

\author{Gareth A. Jones}

\address{School of Mathematical Sciences, University of Southampton, Southampton SO17 1BJ, UK}
\email{G.A.Jones@maths.soton.ac.uk}


%11F06 Modular groups, arithmetic groups
%11G32 Dessins d'enfants, Belyi theory
%14H57 Dessins d'enfants thy
%20B25 Fin auto gps of alg, geom or comb structures
%20H10 Fuchsian groups
%30F10 Compact Riemann surfaces, uniformisation

\maketitle

%%%%%%%%%%%%%%%%%
%%%%%%%%%%%%%%%%%

\begin{abstract}
Greenberg proved that every countable group $A$ is isomorphic to the automorphism group of a Riemann surface, which can be taken to be compact if $A$ is finite. We give a short and explicit algebraic proof of this for finitely generated groups $A$. 
\end{abstract}

%%%%%%%%%%%%%%%%%%%%%%%%%%%%
%%%%%%%%%%%%%%%%%%%%%%%%%%%%

\section{Introduction}\label{intro}

In 1960 Greenberg~\cite{Gre60} proved that every countable group $A$ is isomorphic to the automorphism group of a non-compact Riemann surface, which can be taken to have finite type, that is, to have a finitely generated fundamental group, if $A$ is finite. His proof is quite complicated, using notions of $N$-equivalence and $N$-maximality which he introduced and developed in~\cite{GreCan} (however, see~\cite{All06} for a more elementary geometric proof by Allcock). In 1973 he proved in~\cite[Theorem $6'$]{Gre73} that every finite group is isomorphic to the automorphism group of a compact Riemann surface. (This was also stated without proof in~\cite[Theorem 4]{Gre63}.) His proof of this depends on a delicate construction~\cite[Theorem 4]{Gre73} of maximal Fuchsian groups with a given signature. Here we give a short algebraic proof of his results in the case of finitely generated groups $A$, based on well-known properties of triangle groups and their finite quotient groups (see Theorem~\ref{mainthm}). The author is grateful to Alexander Mednykh for asking whether such a proof might be possible, and to David Singerman for many helpful comments concerning Fuchsian groups. 

%%%%%%%%%%%

\section{The proof}\label{proof}

We will prove the following restricted version of Greenberg's Theorem:

\begin{theo}[Greenberg]\label{mainthm}
Every finitely generated group $A$ is isomorphic to the automorphism group of a Riemann surface, which can be taken to be compact if $A$ is finite.
\end{theo}

\noindent{\sl Proof.} Let $\Delta$ be a hyperbolic triangle group
\[\Delta(l,m,n)=\langle X, Y, Z\mid X^l=Y^m=Z^n=XYZ=1\rangle,\]
so that $l^{-1}+m^{-1}+n^{-1}<1$ and $\Delta$ acts by isometries on the hyperbolic plane $\mathbb H$.
Dirichlet's Theorem on primes in arithmetic progressions  implies that there are infinitely many prime powers $q\equiv -1$ mod~$(k)$ where $k:={\rm lcm}(2l, 2m, 2n)$. For any such $q$ there is a smooth (surface-kernel) epimorphism $\Delta\to G:=PS\negthinspace L_2({\mathbb F}_q)$, so that the images $x, y$ and $z$ of $X, Y$ and $Z$ have orders $l, m$ and $n$ (see~\cite[Corollary C]{Gar15} or~\cite{Macb69}, for example).
These orders divide $(q+1)/2$, so $x, y$ and $z$ act semiregularly in the natural action of $G$ on the projective line ${\mathbb P}^1({\mathbb F}_q)$. Since this action is primitive, the subgroup $H=G_{\infty}$ of $G$ fixing $\infty$ is a maximal subgroup of index $q+1$ in $G$, and hence its inverse image $N$ in $\Delta$ is a maximal subgroup of index $q+1$ in $\Delta$. Since $X, Y$ and $Z$ induce semiregular permutations of orders $l, m$ and $n$ on the cosets of $N$ in $\Delta$, none of their non-identity powers are conjugate to elements of $N$. Thus $N$ has no elliptic elements, so being cocompact (since $\Delta$ is), it is a surface group
\[N=\langle A_i, B_i\; (i=1,\ldots, g) \mid \prod_i[A_i,B_i]=1\rangle\]
of genus $g$ given by the Riemann--Hurwitz formula
\begin{equation}\label{RH}
g=\frac{q+1}{2}\left(1-\frac{1}{l}-\frac{1}{m}-\frac{1}{n}\right)+1.
\end{equation}

Clearly, not every generator $A_i$ or $B_i$ of $N$ can be contained in the core $K$ of $N$ in $\Delta$, since this is the kernel of the action of $\Delta$ on ${\mathbb P}^1({\mathbb F}_q)$ and $N$ acts non-trivially. Without loss of generality (renaming generators if necessary) we may assume that $B_1\not\in K$.

Now $g>(q+1)/84$ by (\ref{RH}), so given any finitely generated group $A$ one can choose $q$ so that $g\ge d$, where $d$ is the rank (minimum number of generators) of $A$. One can then find an epimorphism $\theta:N\to A$ by sending the generators $A_i$ of $\Delta$ to a generating set for $A$, and the generators $B_i$ to the identity.

Let $M=\ker\theta$, so $M$ is normal in $N$ with $N/M\cong A$. Clearly $N_{\Delta}(M)\ge N$, so by the maximality of $N$ in $\Delta$ we must have $N_{\Delta}(M)=N$ or $\Delta$. In the latter case $M$ is normal in $\Delta$ and is therefore contained in $K$, which is impossible since $B_1\in M\setminus K$. Thus $N_{\Delta}(M)=N$.

The argument so far, which applies to {\sl any\/} hyperbolic triple $(l, m, n)$, shows that $A$ is isomorphic to the automorphism group ${\rm Aut}\,{\mathcal M}\cong N_{\Delta}(M)/M=N/M$ of the oriented hypermap $\mathcal M$ of type $(l, m, n)$ corresponding to the subgroup $M$ of $\Delta$. Indeed, if $A$ is finite then so is $|\Delta:M|$, so that the underlying Riemann surface ${\mathcal S}={\mathbb H}/M$ is compact, of genus $|A|(g-1)+1$, and $\mathcal M$ is a dessin d'enfant in Grothendieck's sense~\cite{Gro}, with automorphism group $A$. (See~\cite{JW} for background on hypermaps and dessins.) However, we wish to realise $A$ as the automorphism group of ${\mathcal S}$, rather than $\mathcal M$; certainly ${\rm Aut}\,{\mathcal S}$ contains ${\rm Aut}\,{\mathcal M}$, but it could be larger. Since $N$ has no elliptic elements, neither has $M$, so $M$ acts without fixed points on $\mathbb H$; it follows from this (see~\cite[Theorem~5.9.4]{JS87}, for instance) that ${\rm Aut}\,{\mathcal S}\cong N(M)/M$ where $N(M)$ is the nornaliser of $M$ in ${\rm Aut}\,{\mathbb H}=PS\negthinspace L_2({\mathbb R})$. (Note that $N(M)$ is a Fuchsian group since $M$ is a non-cyclic Fuchsian group, see~\cite[Theorem~5.7.5]{JS87} for instance.) Clearly $N_{\Delta}(M)\le N(M)$, and we need to prove equality here.

In order to do this, let us choose the triple $(l,m,n)$ so that $\Delta$ is maximal (as a Fuchsian group) and non-arithmetic. (By results of Singerman~\cite{Sin72} and Takeuchi~\cite{Tak77} these conditions are satisfied by `most' hyperbolic triples: see Remark~2, following this proof, for specific examples.) Since $\Delta$ is non-arithmetic, a theorem of Margulis~\cite[\S IX.7]{Mar} implies that its commensurator $\overline{\Delta}$ in $PS\negthinspace L_2({\mathbb R})$ is a Fuchsian group. Since $\overline{\Delta}$ contains $\Delta$, the maximality of $\Delta$ implies that $\overline{\Delta}=\Delta$ and hence $\Delta$ is the commensurator of each of its subgroups of finite index, including $N$.

Now let $g\in N(M)$, so $\Delta^g\ge N^g\triangleright M^g=M$ and $N(M)$ contains both $N$ and $N^g$; moreover, it must do so with finite index, since these groups are cocompact, so $N\cap N^g$ has finite index in $N$ and $N^g$. Since $\Delta$ is the commensurator of $N$ it follows that $g\in\Delta$, so $g\in N_{\Delta}(M)$ as required.

%%%%%%%%%%%%%

\section{Remarks}\label{remarks}

\noindent{\bf  1} One cannot regard this as an elementary proof of Greenberg's results (Allcock gives one in~\cite{All06}), since those of Margulis, Singerman and Takeuchi which it uses are far from elementary. Nevertheless, the route from them to the required destination is both short and straightforward.

\smallskip

\noindent{\bf  2} As an example of a triple for which $\Delta$ is non-arithmetic and maximal, one could take $(2,3,n)$ for any prime $n\ge 13$ (or indeed any integer $n>30$). If $n=13$, for instance, we require $q\equiv-1$ mod~$(156)$; the smallest such prime power is the prime $311$, giving genus $g=15$, so that all groups $A$ of rank $d\le 15$ are realised. Taking triples $(2,3,21)$ or $(2,4,9)$ allows smaller primes $q=83$ or $71$, both giving $g=6$. The triple $(4,6,12)$ allows an even smaller prime $q=23$, but leads to a larger genus $g=7$.

\smallskip
\iffalse
\noindent{\bf 2} A slightly weaker lower bound $p\ge\max\{12d-1, 13\}$ can be obtained by using only primes $p\equiv 11$ mod~$(12)$. Thus if $d=2$ one can take $p=23$. Even smaller values of $p$ may be available if $A$ has generators of order $2$ or $3$. For example, every finite simple group $A$ has generators $a$ and $b$ with $b^2=1$ (see~\cite{Kin17} for a rather stronger result by King); if $p=17$ then $N$ has genus $1$ and elliptic periods $2, 2, 17$, so one can realise $A$ by mapping $A_1\mapsto a$ and $Y_1, Y_2\mapsto b$, with $B_1, Z_1\mapsto 1$.

\smallskip

\fi

\noindent{\bf 3} In the above proof $g$ (and hence $d$) has linear growth as $q\to\infty$. Faster growth can be obtained by replacing the natural representation of $G=PS\negthinspace L_2(q)$ with a different primitive representation. For example, if $3<q\equiv\pm 3$ or $\pm 13$ mod~$(40)$ then $G$ has a conjugacy class of maximal subgroups $H\cong A_4$ of  index $q(q^2-1)/24$ (see~\cite[Ch.~XII]{Dic}, for example). If $l, m$ and $n$ are coprime to $6$ then no non-identity powers of $x, y$ or $z$ are conjugate to elements of $H$, so the inverse image $N$ of $H$ in $\Delta$ is a surface group of genus 
\[g=\frac{q(q^2-1)}{48}\left(1-\frac{1}{l}-\frac{1}{m}-\frac{1}{n}\right)+1>\frac{q(q^2-1)}{120}\]
(since $l, m, n\ge 5$), giving cubic growth of $g$ as $q\to\infty$. Again one must choose $(l,m,n)$ so that $\Delta$ is maximal and non-arithmetic: the smallest such example in this case is $(7,11,13)$. One must also choose $q$ so that $G$ has generators $x, y$ and $z$ of orders $l$, $m$ and $n$: for examples, see~\cite{Gar15} or~\cite{Macb69}.

\smallskip

\noindent{\bf 4} There are many other possibilities for $\Delta$, $G$ and $H$ in this proof: one example, giving even faster growth of $g$, is to use the action of the symmetric group $G=S_{\negthinspace\negthinspace p}$, for primes $p\ge 5$, on the $(p-2)!$ cosets of its subgroup $H=AGL_1(p)$ (maximal by the classification of finite simple groups, which includes that of groups of prime degree). Now $H$ is a Frobenius group, that is, non-identity elements of $H$ have at most one fixed point in the natural action of $G$ of degree $p$, so one can ensure that $N$ is a surface group by choosing generators $x, y$ and $z$ for $G$ each with at least two fixed points. One can choose such a triple to generate $G$ as follows. Provided $G_0:=\langle x, y, z\rangle$ is transitive it is primitive since the degree $p$ is prime, so if at least one of $x, y$ and $z$ is a cycle with at least three fixed points then an extension of Jordan's Theorem in~\cite{Jon14} implies that $G_0$ contains the alternating group $A_p$. If, in addition, at least one of  $x, y$ and $z$ is odd then $G_0=S_{\negthinspace\negthinspace p}$, as required. If $l, m$ and $n$ are the orders of $x, y$ and $z$ then the inverse image $N$ of $H$ in $\Delta=\Delta(l,m,n)$ is a surface group of genus
\[g=\frac{(p-2)!}{2}\left(1-\frac{1}{l}-\frac{1}{m}-\frac{1}{n}\right)+1>\frac{(p-2)!}{84},\]
giving super-exponential growth. As an example, if $p=2l-3$ one could take $x=(1, 2, \ldots, l)$ and take $y$ to be the cycle $(4, 3, 2, 1, l+1, l+2, \ldots, p)$ of length $m=l+1$, so that $z=(xy)^{-1}$ is the cycle $(p, p-1, \ldots, 4)$ of length $n=p-3=2l-6$; in this case, since $l, m, n\to\infty$ with $p$, we have $g\sim(p-2)!/2$ as $p\to\infty$. The lists of exceptions in~\cite{Sin72} and~\cite{Tak77} show that here $\Delta(l,m,n)$ is maximal and non-arithmetic for each prime $p\ge 13$; for instance, if $p=13$ then $\Delta=\Delta(8, 9, 10)$, with $g=13250161$ large enough to realise most finitely generated groups of current interest.
%{\color{blue}[Omit this Remark?]}

\smallskip

\noindent{\bf 5} The proof in Section~\ref{proof} is adapted from one in~\cite[Theorem~3(a)]{Jon18} that for any hyperbolic triple $(l,m,n)$ there are $\aleph_0$ non-isomorphic dessins (finite oriented hypermaps) of type $(l,m,n)$ with a given finite automorphism group $A$. (See~\cite{Hid} for related results by Hidalgo on realising groups as automorphism groups of dessins.) The Riemann surfaces $\mathcal S$ underlying these dessins have automorphism group containing $A$. It would be interesting to determine whether they can be chosen so that ${\rm Aut}\,{\mathcal S}=A$ in those cases where the corresponding triangle group is arithmetic or non-maximal.

\smallskip

\noindent{\bf 6} In~\cite{Jon18} it is also shown that for many hyperbolic triples (including all of non-cocompact type and many of cocompact type), every countable group can be realised as the automorphism group of $2^{\aleph_0}$ non-isomorphic oriented hypermaps of that type. It would be interesting to try to deduce Greenberg's more general theorem for all countable groups $A$ from this type of argument. The main difficulty is that one needs $N$ to have infinite rank, and hence to have infinite index in $\Delta$, so that commensurability cannot be used as in the last paragraph of Section~\ref{proof}. (See also~\cite{GTW, SS} for similar applications of commensurability, and for examples of how arithmetic and non-arithmetic triangle groups behave differently.)
%{\color{blue}[But see \S\ref{Resid}, if it exists.]}

\smallskip

\noindent{\bf 7} In all the above variations of this proof, if $A$ is finite then since the subgroups $M$ have finite index in triangle groups $\Delta$, Bely\u\i's Theorem~\cite{Bel}, as reinterpreted by Grothendieck~\cite{Gro}, implies that the Riemann surfaces ${\mathcal S}={\mathbb H}/M$ are defined, as projective algebraic curves, over algebraic number fields. (See~\cite[Ch.~1]{JW} for background on Bely\u\i's Theorem.) What can be said about these fields? For example, can every finite group $A$ be realised as the automorphism group of a curve (or dessin) defined over $\mathbb Q$? By Cayley's Theorem, $A$ is contained in such a group: the standard generating triples for $S_{\negthinspace n}$, consisting of cycles of lengths $2$, $n-1$ and $n$, are mutually conjugate, so they correspond to a unique regular dessin $\mathcal D$ with automorphism group $S_{\negthinspace n}$; the absolute Galois group ${\rm Gal}\,\overline{\mathbb Q}/{\mathbb Q}$ (where $\overline{\mathbb Q}$ is the field of algebraic numbers) preserves the automorphism group and passport (triple of cycle-structures of generators) of any dessin~\cite{JSt}, so it preserves $\mathcal D$ and hence $\mathcal D$ is defined over $\mathbb Q$. (See~\cite[\S4.4, \S7.4.1, \S8.3.1]{Ser} for a Galois-theoretic interpretation by Serre of this example of rigidity.)


\begin{thebibliography}{99}

\bibitem{All06}
D.~Allcock,
Hyperbolic surfaces with prescribed infinite symmetry groups,
{\sl Proc.~Amer.~Math.~Soc.} 134 (2006), 3057--3057.

\bibitem{Bel}
G.~V.~Bely\u\i,
On Galois extensions of a maximal cyclotomic field,
{\sl Izv.~Akad.~Nauk SSSR Ser.~Mat.} 43 (1979), 267--276, 479.

\bibitem{Dic}
L.~E.~Dickson,
{\sl Linear Groups},
Dover, New York, 1958.

\bibitem{Gar15}
S.~Garion,
On Beauville structures for $PS\negthinspace L_2(q)$,
{\sl J.~Group Theory} 18 (2015), 981--1019.

\bibitem{GTW}
E.~Girondo, D.~Torres-Teigell and J.~Wolfart,
Shimura curves with many uniform dessins,
{\sl Math.~Z.} 271 (2012), 757--779.

\bibitem {Gre60}
L.~Greenberg,
Conformal transformations of Riemann surfaces,
{\sl Amer.~J.~Math.} 82 (1960), 749--760.

\bibitem{GreCan}
L.~Greenberg,
Discrete groups of motions,
{\sl Canad.~J.~Math.} 12 (1960), 415--426.

%\iffalse
\bibitem{Gre63}
L.~Greenberg,
Maximal Fuchsian groups,
{\sl Bull.~Amer.~Math.~Soc.} 69 (1963), 569--573.
%\fi

\bibitem {Gre73} 
L.~Greenberg,
Maximal groups and signatures,
in {\sl Discontinuous Groups and Riemann Surfaces, Proc.~Conf. Maryland 1973},
ed.~L.~Greenberg, Princeton University Press, Princeton NJ, 1974, 207--226.

\bibitem{Gro}
A. Grothendieck,
Esquisse d'un Programme (1984),
in {\sl Geometric Galois Actions 1},
%: Around Grothendieck's Esquisse d'un Programme
ed.~L. Schneps and P. Lochak, London Math.~Soc.~Lecture Note Ser.~242,
Cambridge University Press, Cambridge, 1997, 5--47.

\bibitem{Hid}
R.~Hidalgo,
Automorphism groups of dessins d'enfants,
{\sl Arch.~Math.~(Basel)} 112 (2019), 13--18.

\bibitem{Jon18}
G.~A.~Jones,
Realisation of groups as automorphism groups in categories,
arXiv:1807.00547v3 [math.GR].

\bibitem{Jon14}
G.~A.~Jones,
Primitive permutation groups containing a cycle,
{\sl Bull~Aust.~Math.~Soc.} 89  (2014), 159--165.

\bibitem{JS87}
G.~A.~Jones and D.~Singerman,
{\sl Complex Function Theory},
Cambridge University Press, Cambridge, 1987.

\bibitem{JSt}
G.~A.~Jones and M.~Streit,
Galois groups, monodromy groups and cartographic groups,
in {\sl Geometric Galois Actions 2},
%: The Inverse Galois Problem, Moduli Spaces and Mapping Class Groups},
ed.~P.~Lochak and L.~Schneps, London Math.~Soc.~Lecture Note Ser.~243,
Cambridge University Press, Cambridge, 25--65.

\bibitem{JW}
G.~A.~Jones and J.~Wolfart,
{\sl Dessins d'Enfants on Riemann Surfaces},
Springer, Cham, 2016.

%\bibitem{Kin17}
%C.~S.~H.~King,
%Generation of finite simple groups by an involution and an element of prime order,
%{\sl J.~Algebra} 478 (2017), 153--173.

\bibitem{Macb69}
A.~M.~Macbeath,
Generators of the linear fractional groups, in {\sl Number Theory (Houston 1967)},
ed.~W.~J.~Leveque and E.~G.~Straus, Proc.~Sympos.~Pure Math.~12,
Amer.~Math.~Soc., Providence RI, 1969, 14--32.

\bibitem{Mar}
G.~Margulis,
{\sl Discrete Subgroups of Semisimple Lie Groups},
Springer-Verlag, Berlin, 1991.

%\bibitem{Mat} B.~H.~Matzat, {\sl Konstruktive Galoistheorie}, Springer Lecture Notes in Math. 1284, Springer, Berlin, 1986.

\bibitem{Ser}
J-P.~Serre,
{\sl Topics in Galois Theory},
Jones and Bartlett, Boston Mass., 1992

%\bibitem{Sin70} D.~Singerman, Subgroups of Fuchsian groups and finite permutation groups,
%{\sl Bull.~London Math.~Soc.} 2 (1970), 319--323.

\bibitem{Sin72} 
 D.~Singerman,
 Finitely maximal Fuchsian groups,
 {\sl J.~London Math.~Soc.~(2)} 6 (1972), 29--38.
 
 \bibitem{SS} 
 D.~Singerman and R.~I.~Syddall,
 The Riemann surface of a uniform dessin,
 {\sl Beitr\"age Algebra Geom.} 44 (2003), 413--430.
 
 \bibitem{Tak77}
 K.~Takeuchi,
 Arithmetic triangle groups,
 {\sl J.~Math.~Soc.~Japan} 29 (1977), 91--106.

\end{thebibliography}
\end{document}